\documentclass[12pt]{amsart}
\usepackage{amssymb}

\renewcommand{\baselinestretch}{1.5}
\theoremstyle{plain}
\newtheorem{theorem}{Theorem}

\theoremstyle{definition}
\newtheorem{definition}[theorem]{Definition}

\newcommand\tr{\operatorname{tr}}
\newcommand\rng{\operatorname{rng}}

\begin{document}
\vglue -20pt

\centerline{\large{\textbf{PRESERVING THE MEASURE OF
COMPATIBILITY}}}

\centerline{\large{\textbf{BETWEEN QUANTUM STATES}}}

\vskip 20pt \centerline{LAJOS MOLN\'AR}

\vskip  0pt \centerline{Institute of Mathematics and Informatics}

\vskip -5pt \centerline{University of Debrecen}

\vskip -5pt \centerline{4010 Debrecen, P.O.Box 12}

\vskip -5pt \centerline{Hungary}

\vskip -5pt \centerline{e-mail: \texttt{molnarl@math.klte.hu}}

\vskip  5pt \centerline{and}

\vskip 10pt \centerline{WERNER TIMMERMANNN}

\vskip  0pt \centerline{Institut f\"ur Analysis}

\vskip -5pt \centerline{Technische Universit\"at Dresden}

\vskip -5pt \centerline{D-01062 Dresden}

\vskip -5pt \centerline{Germany}

\vskip -5pt \centerline{e-mail: \texttt{timmerma@math.tu-dresden.de}}

\vskip 20pt \centerline{Running title:}

\centerline{\scshape{Preserving the measure of compatibility}}

\pagestyle{myheadings}
\markboth{\textsc{\SMALL Preserving the measure of compatibility}}
{\textsc{\SMALL Preserving the measure of compatibility}}

\normalsize \vskip 20pt
\centerline{\textsc{Abstract}}
In this paper after defining the abstract concept of
compatibility-like functions on quantum states, we prove that
every bijective transformation on the set of all states which
preserves such a function is implemented by an either unitary or
antiunitary operator.

\newpage
\renewcommand{\baselinestretch}{1.5}
\normalsize

In the last couple of years several communications have appeared
in connection with the following problem raised by R. Peierls
\cite{Pei1}: when different density matrices can characterize the
knowledge available to different people about one and the same
physical system? The first answer given by Peierls in
\cite{Pei1,Pei2} was that the density matrices under consideration
must commute and their product must be nonzero. However, C. Fuchs
\cite{Fuc} (also see \cite{Mer}) gave an example which made
Peierls' first condition questionable. After that several attempts
have been made to find the proper solution of the problem (e.g.,
\cite{BFM, Mer}). All those attempts operate with the concept of
compatibility of density matrices. According to them, we say that
a collection of density matrices is compatible if the supports of
the matrices under consideration (i.e., the orthogonal complements
of their null spaces) have nontrivial intersection. So, it is just
an easy task to determine whether a pair of density matrices is
compatible or not. Having this in mind, it is now a natural
problem to give sense to the following question: if a pair of
density matrices is compatible, then "how much" compatible they
are. In other words, we arrive at the problem of measuring the
compatibility. One possibility to define such a measure was
described in \cite {PoBK}. Namely, in some analogy with the
fidelity, C. Poulin and R. Blume-Kohout defined a compatibility
function \cite[Definition 1]{PoBK} which fulfils certain natural
physical requirements and they proved some important properties.

In our recent paper \cite{Mol2} we have determined the structure
of the bijective transformations on the set of all density
operators which preserve the fidelity. This result is in close
relation with Wigner's theorem on symmetry transformations. In
fact, it can be considered as a Wigner-type result for the set of
all mixed states (recall that Wigner's original result concerns
the pure states). In \cite{Mol2} we proved that the
transformations in question are all implemented by unitary or
antiunitary operators on the underlying Hilbert space. In view of
this result and the analogy between the fidelity and the measure
of compatibility defined by Poulin and Blume-Kohout, it is a
natural problem to determine the structure of the bijective
transformations of the set of all density matrices which preserve
the compatibility function. We shall see below that the solution
of this problem is the same as the one concerning fidelity. This
is the content of the present paper.

Let us begin with the notation. Let $H$ be a (complex, not
necessarily finite dimensional) Hilbert space. If not stated
otherwise, all operators on $H$ are meant to be bounded and
linear. The expression $\rng A$ denotes the range of the operator
$A$. If $A$ is positive, $A^{1/2}$ stands for its unique positive
square root.

Denote by $S(H)$ the set of all states (or, in another
terminology, density operators) on $H$ i.e., the positive
trace-class operators on $H$ with trace 1. The set $S(H)$ is a
convex subset of the space of all self-adjoint operators on $H$
and its extreme points (which are exactly the rank-one
projections) are called pure states.

Now, instead of using the concept due to Poulin and Blume-Kohout,
we define the abstract concept of compatibility-like functions
which extends \cite[Definition 1]{PoBK} to obtain a result of
higher generality.

\begin{definition}\label{D:compat}
Let $C:S(H)\times S(H)\to [0,1]$ be a function such that for any
pair $A,B \in S(H)$ of states we have
\begin{itemize}
\item[(i)] $C(A,B)=0$ if and only if $\rng A^{1/2}\cap \rng
B^{1/2}=\{ 0\}$,
\item[(ii)] $C(A,B)=C(B,A)$,
\item[(iii)] if $P$ is a pure state, then
\[
C(A,P)^2=\sup \{ \lambda \in [0,1] \, : \, \lambda P\leq A\}.
\]
\end{itemize}
We say that $C$ is a compatibility-like function on the set of all
states on the Hilbert space $H$.
\end{definition}

Several remarks should be made concerning the above definition.
First, we emphasize that our definition is formulated for both
finite and infinite dimensional Hilbert spaces (in
\cite[Definition 1]{PoBK} only finite dimensional spaces were
considered). Concerning the correctness of the definition we note
the following. The quantity on the right hand side of the equality
in (iii) also appears in relation with effects. A self-adjoint
operator $T$ on $H$ with the property $0\leq T\leq I$ ($I$ is the
identity operator) is called an effect. The effects are well-known
to play important role in the quantum theory of measurement (e.g.,
\cite{BusLahMit}). Now, it is clear that every state on $H$ as a
linear operator can also be viewed as an effect. If $T$ is an
effect, $\varphi$ is a unit vector in $H$ and $P_\varphi$ denotes
the rank-one projection onto the subspace generated by $\varphi$,
then the quantity
\[
\lambda(T, P_\varphi)= \sup \{ \lambda \in [0,1] \, : \, \lambda
P_\varphi \leq T\}
\]
is called the strength of $T$ along the ray represented by
$\varphi$. This concept was introduced by Busch and Gudder in
\cite{BuGu}. It was proved in \cite[Theorem 3]{BuGu} that
$\lambda(T, P_\varphi)=0$ if and only is $\varphi \notin \rng
T^{1/2}$ which is equivalent to $\rng T^{1/2} \cap \rng P_\varphi^{1/2}
=\{0\}$. This means that there is no contradiction between the
conditions (i) and (iii).

Observe that if $H$ is finite dimensional, then the ranges of a
positive operator and its square root are the same and they are
automatically closed. Therefore, in the finite dimensional case we
have
\[
\rng A^{1/2}=\rng A={\overline{\rng A}}=(\ker A)^\perp
\]
and hence (i) says that $C(A,B)>0$ if and only if $A,B$ are
compatible in the sense mentioned in the introduction. The meaning
of (ii) is clear. Now, what about (iii)? One might think that this
condition is quite restrictive and probably has no physical
meaning for states. But it can
be shown that the compatibility function defined by Poulin and
Blume-Kohout satisfies (iii) (see either \cite[Definition 1]{PoBK}
itself or \cite[Theorem 3]{PoBK}) as well as (i) and (ii). So, to
sum up, our definition is a generalization of the one given by
Poulin and Blume-Kohout and hence it certainly has sense at least
from the mathematical point of view.

The reason why we assume (iii) is that there is a
nice formula to compute $C(A,P)$. Namely, by \cite[Theorem
4]{BuGu} for every unit vector $\varphi\in H$ we have
\begin{equation}\label{E:comp4}
C(A,P_\varphi)^2= \lambda(A, P_\varphi)=
\left\{%
\begin{array}{ll}
    \|A^{-1/2}\varphi\|^{-2}, & \hbox{{\text{if }} $\varphi \in \rng (A^{1/2})$;} \\
    0, & \hbox{{\text{else}}.} \\
\end{array}%
\right.
\end{equation}
(Here $A^{-1/2}$ denotes the inverse of $A^{1/2}$ on $\rng
A^{1/2}$.) The proof of our result is based on this
correspondence.

We further note that it would be another natural assumption to
suppose that $C$ is invariant under unitary-antiunitary
transformations. But, as we do not need it in our proof, we do not
assume it.

To conclude our remarks, we show a natural example for a
compatibility-like function which might also justify our
definition. So, for any pair $A,B \in S(H)$ define
\[
\begin{aligned}
C(A,B) = \sup \{ & \sum_n \sqrt{\lambda_n \mu_n} \, :\,
\lambda_n,\mu_n \in [0,1], \, \sum_n \lambda_n=\sum_n \mu_n=1
\text{ and }
\\ &
\exists \text{ pure states } Q_n \text{ with } \sum_n \lambda_n
Q_n=A, \, \sum_n \mu_n Q_n=B\}.
\end{aligned}
\]
It is easy to verify that this function has the properties
(i)-(iii).  In fact, in accordance with the discussions in
\cite{BFM} and \cite{Mer}, we believe that this compatibility-like
function represents the most natural way of defining a measure of
compatibility between quantum states.

Now, our result reads as follows.

\begin{theorem}
Let $H$ be a Hilbert space and let $C$ be a compatibility-like
function on $S(H)$. Let $\phi: S(H) \to S(H)$ be a bijective
function which preserves $C$, that is, assume that
\[
C(\phi(A),\phi(B))=C(A,B) \qquad (A,B\in S(H)).
\]
Then there exists an either unitary or antiunitary operator $U$ on
$H$ such that $\phi$ is of the form
\[
\phi(A)=UAU^* \qquad (A\in S(H)).
\]
\end{theorem}

\begin{proof}
Clearly, we can assume that $\dim H\geq 2$. For temporary use, we
say that the states $D,A$ are compatible (resp. incompatible) if
$C(D,A)>0$ (resp. $C(D,A)=0$). It is useful to introduce the
following notation. If $\mathcal M$ is a subset of $S(H)$, then
denote
\[
{\mathcal M}^{ic}=\{ D\in S(H) \, :\, C(D, A)=0 \text{ for all }
A\in \mathcal M\}.
\]
By condition (i) in Definition~\ref{D:compat}, $C(D,A)=0$ means
that the subspaces $\rng D^{1/2}$, $\rng A^{1/2}$ of $H$ have
trivial intersection. It can be easily verified that the operator
$A\in S(H)$ has rank one (which means that $A$ is a pure state) if
and only if
\[
{({\{A\}}^{ic})}^{ic}=\{ A\}.
\]
Since $\phi$ preserves the compatibility in both directions, it
follows from this characterization that $\phi$ preserves the pure
states in both directions. This means that $A\in S(H)$ is a pure
state if and only if so is $\phi(A)$.

Next we assert that $\phi$ maps independent pure states to
independent ones. Here, a set of $n$ pure states (rank-one
projections) is called independent if their ranges generate an
$n$-dimensional subspace of $H$. To prove the assertion we use
induction. The statement is obvious if the set has only one
element. Let $\{P_1,\ldots P_n, P_{n+1}\}$ be a set of $n+1$ pure
states such that the subset $\{P_1, \ldots, P_n\}$ is independent.
It is easy to see that $\{P_1,\ldots P_n, P_{n+1}\}$ is dependent
if and only if for any $A\in S(H)$ with $C(A,P_1)>0, \ldots,
C(A,P_n)>0$ we have that $C(A,P_{n+1})>0$ holds too. Indeed, this
follows from the fact that for any pure state $P$ we have
$C(A,P)>0$ if and only if the range of $P$ is included in the
range of $A^{1/2}$ (see the remarks after the
Definition~\ref{D:compat}). Using the above description of
dependence, it is now clear that assuming $\phi$ maps independent
sets of $n$ pure states to sets of the same kind, we have the same
property of $\phi$ for $n+1$ in the place of $n$. Since
$\phi^{-1}$ has the same properties as $\phi$, we deduce that
$\phi$ preserves the independence of the sets of pure states in
both directions.

It is easy to see that an operator $A\in S(H)$ has rank $n$ if and
only if there exists an independent set of $n$ pure states such
that $C(A,P)>0$ for every element of that set, but there does not
exist a set of $n+1$ elements having the same property. This gives
us that $\phi$ preserves the rank.

Now we prove that $\phi$ preserves the transition probability
between pure states. Recall that for any pair $P,Q$ of pure
states, the transition probability between them is $\tr PQ$, where
$\tr$ is the usual trace-functional.  To verify the mentioned
preserver property of $\phi$, first let $P,Q$ be rank-one
projections with orthogonal ranges. Define
\[
A=\lambda P+\mu Q,
\]
where $\lambda, \mu$ are fixed and satisfy $0<\lambda <\mu <1,
\lambda+\mu=1$. Clearly, $A$ acts on the 2-dimensional subspace
$H_A$ of $H$ generated by the ranges of $P,Q$. (Here, the phrase
that $A$ acts on $H_A$ means that $(\ker A)^\perp =
{\overline{\rng A}}=\rng A=H_A$.) We assert that $\phi(A)$ acts on
the subspace generated by the ranges of the independent pure
states $\phi(P), \phi(Q)$. Indeed, $\phi(A)$ has rank 2 and taking
into account that
\[
C(\phi(A),\phi(P))=C(A,P)>0, \quad C(\phi(A),\phi(Q))=C(A,Q)>0,
\]
we see that the ranges of $\phi(P), \phi(Q)$ are included in $\rng
\phi(A)^{1/2}= \rng \phi(A)$. This clearly implies our assertion.
In what follows we restrict the considerations onto those
2-dimensional subspaces, that is, to the ranges of $A$ and
$\phi(A)$, respectively. By property (iii) in the definition of
compatibility-like functions, we see that
\[
\lambda \leq C(A,R)^2 \leq \mu
\]
holds for every rank one projection $R$ on the range of $A$. As
$\phi$ preserves $C$, we have
\[
\lambda \leq C(\phi(A), \phi(R))^2 \leq \mu
\]
for every rank one projection $\phi(R)$ on the range of $\phi(A)$.
Moreover, we have
\[
C(\phi(A),\phi(P))^2=C(A,P)^2=\lambda, \quad
C(\phi(A),\phi(Q))^2=C(A,Q)^2=\mu.
\]
Now we refer to a result in \cite{Mol1}. Namely, Lemma 3 given
there states that if $T$ is an effect and $0<\epsilon< \delta\leq
1$ are scalars such that $\epsilon I\leq T\leq \delta I$ and we
have unit vectors $\varphi, \psi\in H$ such that
$\lambda(T,P_\varphi)=\epsilon$ and $\lambda(T,P_\psi)=\delta$,
then $\varphi,\psi$ are eigenvectors of $T$ and the corresponding
eigenvalues are $\epsilon, \delta$, respectively. Using this
result and the correspondence between compatibility-like functions
and the strength, we obtain that the range of $\phi(P)$ is the
eigensubspace of $\phi(A)$ corresponding to the eigenvalue
$\lambda$ and the range of $\phi(Q)$ is the eigensubspace of
$\phi(A)$ corresponding to the eigenvalue $\mu$. Therefore, we
have
\begin{equation}\label{E:comp1}
\phi(A)=\lambda \phi(P)+\mu \phi(Q).
\end{equation}

Now let $P, R$ be arbitrary rank-one projections. Pick a rank-one
projection $Q$ which is orthogonal to $P$ such that the subspace
generated by the ranges of $P$ and $Q$ includes the range of $R$.
Let $\lambda,\mu$ and $A$ be as above. It is easy to check that by
the formula \eqref{E:comp4} we have
\begin{equation}\label{E:comp2}
\begin{gathered}
C^2(A, R)=
\frac{1}{\frac{1}{\lambda}\tr PR+\frac{1}{\mu}\tr QR}=\\
\frac{\lambda\mu}{{\mu}\tr PR+{\lambda}\tr QR}=
\frac{\lambda\mu}{{\mu}\tr PR+{\lambda}(1-\tr PR)}=\\
\frac{\lambda\mu}{{(\mu-\lambda)}\tr PR+{\lambda}}.
\end{gathered}
\end{equation}
As the spectral resolution of $\phi(A)$ is \eqref{E:comp1}, we
similarly have
\begin{equation}\label{E:comp3}
C^2(\phi(A), \phi(R))=
\frac{\lambda\mu}{{(\mu-\lambda)}\tr \phi(P)\phi(R)+{\lambda}}.
\end{equation}
Since $\phi$ preserves $C$, it follows from \eqref{E:comp2} and
\eqref{E:comp3} that
\[
\tr \phi(P)\phi(R)=\tr PR,
\]
which means that $\phi$ preserves the transition probability
between pure states. It follows from Wigner's theorem that $\phi$,
when restricted onto the set of all pure states, is of the form
\[
\phi(P)=UPU^*
\]
for some unitary or antiunitary operator $U$ on $H$.

It remains to show that the above formula holds for every state as
well. The proof goes as follows. Let $A\in S(H)$. For every
rank-one projection $P$ we compute
\[
\begin{gathered}
\lambda(UAU^*, P)= \lambda(A, U^* PU)= C(A, U^* PU)^2=\\
C(\phi(A), \phi(U^* PU))^2= C(\phi(A), P)^2=\lambda(\phi(A), P).
\end{gathered}
\]
Now we refer to \cite[Corollary 1]{BuGu} which sates that if the
strengths of two effects are the same along every ray, then the
effects in question are equal. This gives us that
\[
\phi(A)=UAU^* \qquad (A\in S(H))
\]
and the proof is complete.
\end{proof}

\section*{Acknowledgements}
This paper was written when the first author held a Humboldt
Research Fellowship. He is very grateful to the Alexander von
Humboldt Foundation for providing ideal conditions for research
and to the staff of the Institute of Analysis, TU, Dresden
(Germany) for the kind hospitality.
         The first author also acknowledges support from
         the Hungarian National Foundation for Scientific Research
         (OTKA), Grant No. T030082, T031995, and from
         the Ministry of Education, Hungary, Grant
         No. FKFP 0349/2000

\newpage

\end{document}